\documentclass[12pt]{article}
\usepackage{amsmath,amssymb,amsthm,eucal,upref,bm}%

\begin{document}

\makeatletter

\newcommand{\E}{\mathrm{e}\kern0.2pt}
\newcommand{\D}{\mathrm{d}\kern0.2pt}
\newcommand{\RR}{\mathrm{I\kern-0.20emR}}

{\renewcommand{\thefootnote}{} 
\renewcommand{\refname}{{\large Rerferences}}

\title{\bf V.\,A. Steklov and the Problem\\ of Sharp (Exact) Constants\\ in
Inequalities of Mathematical Physics}

\author{A.\,I. Nazarov, N.\,G. Kuznetsov and S.\,V. Poborchi$^*$}

\date{}

\maketitle

\setcounter{equation}{0}

\footnotetext {$^*$Alexander I. Nazarov is a leading researcher at the Laboratory of
Mathematical Physics, St.\,Petersburg Department of the Steklov Mathematical
Institute, and a professor at the Department of Mathematical Physics, Faculty of
Mathematics and Mechanics, St.\,Petersburg State University. His work was supported
by the RFBR grant 11-01-00825 and by the St.\,Petersburg University grant
6.38.670.2013. 

Nikolay G. Kuznetsov heads the Laboratory for Mathematical Modelling of Wave
Phenomena at the Institute for Problems in Mechanical Engineering, Russian Academy
of Sciences, St.\,Petersburg. 

Sergey V. Poborchi is a professor at the Department of Parallel Algorithms, Faculty
of Mathematics and Mechanics, St.\,Petersburg State University. His work was
supported by the RFBR grant HK-11-01-00667/13. 
}
}

The 150th anniversary of the birth of the outstanding Russian mathematician Vladimir
Andreevich Steklov falls on January 9, 2014. All over the world, every active
researcher in any area of mathematics knows this name (at least the last one).
Indeed, widely known mathematical institutes of the Russian Academy of Sciences in
Moscow and St.\,Petersburg (there are several others scattered over the vast Eastern
part of Russian Federation) are named after Steklov, thus commemorating the fact
that he was the founding father of their predecessor --- the Physico-Mathematical
Institute established in 1921 in starving Petrograd (the Civil War still lasted in
some corners of what would become the USSR next year). Steklov was the first
director of the institute until his untimely death on May 30, 1926. In this paper
(it is a part of survey to be published in January 2014), we describe advances in
one of several areas influenced by Steklov's work.\medskip


In 1896, A.\,M.~Lyapunov established that the trigonometric Fourier coefficients of a
bounded function Riemann integrable on $(-\pi, \pi) $ satisfy the closedness
equation. He presented this result at a session of the Kharkov Mathematical Society,
but left it unpublished. (At that time, Lyapunov's main occupation was to
investigate properties of the double layer potential and the Dirichlet problem; see
\cite{Lyap} for the full-length presentation of his results in this area.) The same
year, Steklov had taken up studies of the closedness equation initiated by his
teacher; Steklov's extensive work on this topic lasted for 30 years until his death.
For this reason A.~Kneser \cite{Kneser} referred to this equation as ``Steklov's
favorite formula''. It should be mentioned that the term {\it closedness equation}
was indeed introduced by Steklov for general orthonormal systems, but only in 1910
(see brief announcements \cite{S3} and the full-length paper \cite{S4}).

The same year (1896), Steklov \cite{S1} proved that the following inequality 
\begin{equation} 
\int_0^l u^2 (x) \, \D x \leq \left( \frac{l}{\pi}
\right)^2 \int_0^l [u'(x)]^2 \, \D x 
\label{Pob_1}
\end{equation}
holds for all functions which are continuously differentiable on $[0, l]$ and have
zero mean there. For this purpose he used the closedness equation for the Fourier
coefficients of $u$ (the corresponding system is $\{ \cos \, (k \pi x / l)
\}_{k=0}^\infty$ normalised on $[0, l]$). Inequality (\ref{Pob_1}) was among
earliest inequalities with sharp constant that appeared in mathematical physics. It
was then applied for justifying the Fourier method for initial-boundary value
problems for the heat equation in two dimensions with variable coefficients
independent of time. (Later, Steklov justified the Fourier method for the wave
equation as well.) The fact that the constant in (\ref{Pob_1}) is sharp was
emphasized by Steklov in \cite{S2}, where he gave another proof of this inequality
(see pp.~294--296). There is another result proved in \cite{S2} (see pp.~292--294);
it says that (\ref{Pob_1}) is true for continuously differentiable functions
vanishing at the interval's end-points, and again the constant is sharp. In his
monograph \cite{MonPob}, Steklov presented inequality (\ref{Pob_1}) and another one
slightly generalizing it.

Next year (1897), Steklov published the article \cite{S5}, in which the following
analogue of inequality (\ref{Pob_1}) was proved:
\begin{equation}
\int_\Omega u^2 \, \D x \leq C \int_\Omega |\nabla u|^2 \, \D x . \label{Pob_1'}
\end{equation}
Here $\nabla$ stands for the gradient operator and the integral on the right-hand
side is called the Dirichlet integral. The assumptions made by Steklov are as
follows: $\Omega$ is a bounded three-dimensional domain whose boundary is piecewise
smooth and $u$ is a real ${\cal C}^1$-function on $\bar \Omega$ vanishing on $\partial \Omega$.
Again, inequality (\ref{Pob_1'}) was obtained by Steklov with the sharp constant
equal to $\lambda_1^{-1}$, where $\lambda_1$ is the smallest eigenvalue of the
Dirichlet Laplacian in $\Omega$. In the early 1890s, H. Poincar\'e \cite{Poin1} and
\cite{Poin2} obtained (\ref{Pob_1'}) using different assumptions, namely, $u$ has
zero mean over $\Omega$ which is a union of a finite number of smooth convex two- and
three-dimensional domains, respectively. In the latter case, the sharp constant in
(\ref{Pob_1'}) is again $\lambda_1^{-1}$, but $\lambda_1$ is the smallest positive
eigenvalue of the Neumann Laplacian in $\Omega$.

The problem of finding and estimating sharp constants in inequalities attracted much
attention from those who work in theory of functions and mathematical physics (see,
for example, the classical monographs \cite{HLP} and \cite{PS}). It is worth
mentioning that in \cite{HLP}, sect.~7.7, inequality (\ref{Pob_1}) is proved under
either type of conditions proposed by Steklov, but the authors speak about {\it
Wirtinger's inequality}\/ in the above mentioned section with the reference to
\cite{Bla}, p.~105. This confirms the Arnold Principle: ``If a notion bears a
personal name, then this name is not the name of the discoverer'' (see
\cite{Arnold}).

More than thirty years ago, the role of sharp constants was emphasized in the book
\cite{Mih1} by S.\,G.~Mikhlin (he graduated from the Leningrad University a few
years after Steklov's death; see his recollections of student years \cite{Mih2}).
Let us quote the review \cite{Peet} of the German version of this book. 
\begin{quote}
[This book] is devoted to appraising the (best) constants --- exact results or
explicit (numerical) estimates --- in various inequalities arising in
``analysis'' (=PDE). [\dots] This is a most original work, a bold attack in a
direction where still very little is known.
\end{quote} 
Our aim is to outline main achievements in this area, but we restrict ourselves to
the direct generalizations of (\ref{Pob_1}) and (\ref{Pob_1'}), that is, to
inequalities of the following form: 
\begin{equation}
\|u\|_{L^q (\Omega)} \leq C \, \|\nabla u\|_{L^p(\Omega)} . \label{q-p}
\end{equation}
Here $\Omega$ is a domain in $\RR^n$, $n \geq 1$, whereas $p, q \geq 1$ satisfy the
following restrictions: 
\begin{equation*}
\aligned
&q\le p^*=\frac{np}{n-p},&\quad\mbox{if}&\quad 1\le p<n;\\
&q<\infty,&\quad\mbox{if}&\quad p=n>1;\\
&q\le\infty,&\quad\mbox{if}&\quad p>n \quad\mbox{or}\quad n=1.
\endaligned
\end{equation*}
It is assumed that $u$ belongs to $L^{1,p} (\Omega)$, that is, $u \in L^p_{loc} (\Omega)$,
whereas $\nabla u \in L^p (\Omega)$.

Weighted inequalities --- the Hardy inequality and its generalizations such as the
Hardy--Sobolev inequality, the Maz'ya inequality, the Caffarelly--Kohn--Nirenberg
inequality --- will not be considered here. Inequalities involving derivatives of
higher order (they received much attention during the past few years) are also out
of our scope.\medskip

If $u$ vanishes on $\partial \Omega$ (this is understood as follows: $u$ can be
approximated in the norm $\|\nabla u\|_{L^p (\Omega)}$ by smooth functions having compact
support in $\Omega$), then (\ref{q-p}) is true with some positive constant $C$ for any
domain of finite volume\footnote {This condition is not sharp. In the recent papers \cite{JMV} and
\cite{JMV1}, the necessary and sufficient condition for the validity of (\ref{q-p})
is given in the case when $p=q$.} (for an arbitrary domain in the critical case $p<n$,
$q=p^*$). For these functions, inequality (\ref{q-p}) often appears under various
names for different values of $p$ and $q$. In particular, it is referred to as:
\begin{itemize}
 \item the {\it Steklov inequality} when $p=q=2$;

 \item the {\it Friedrichs inequality} when $p=q$;

 \item the {\it Sobolev inequality} when $p<n$, $q=p^*$.
\end{itemize}
Notice that a slightly different inequality was obtained by K.-O.
Friedrichs \cite{Frie} under the assumption that $\Omega \subset \RR^2$. Namely, his
inequality is as follows:
\begin{equation}
\int_\Omega u^2 \, \D x \leq C \left[ \int_\Omega |\nabla u|^2 \, \D x + \int_{\partial \Omega} u^2
\, \D S \right] , \label{Fried}
\end{equation}
where $\D S$ denotes the element of area of $\partial \Omega$. Generally speaking,
(\ref{Fried}) holds for all bounded domains in $\RR^n$, for which the divergence
theorem is true (see \cite{Maz}, p.~24). Furthermore, the Sobolev inequality was
proved by S.\,L. Sobolev himself only for $p>1$, whereas E. Gagliardo proved it for
$p=1$ (see \cite{Sob} and \cite{Ga1}, respectively).

Inequality (\ref{q-p}) for functions $u$ with zero mean value over $\Omega$ is equivalent
to the following one valid for all $u\in L^{1,p} (\Omega)$: 
\begin{equation}
\| u - \langle u \rangle \|_{L^q (\Omega)} \leq C \|\nabla u \|_{L^p (\Omega)} , \qquad \langle u
\rangle = \frac{\int_\Omega u \, \D x}{{\rm meas}_n \, \Omega} \, . \label{Poi}
\end{equation}
Here the $n$-dimensional measure of $\Omega$ stands in the denominator. Moreover, some
requirements must be imposed on $\Omega$ for the validity of (\ref{Poi}). Indeed, as
early as 1933 O.~Nikod\'ym \cite{Nik} (see also \cite{Maz}, p.~7) constructed a
bounded two-dimensional domain $\Omega$ and a function with the finite Dirichlet integral
over $\Omega$ such that inequality (\ref{Poi}) is not true for $p=q=2$. Another example
of a domain with this property is given in \cite{CH}, ch.~7, sect.~8.2 (see also
\cite{Maz}, sect.~6.10.3). On the other hand, if $p=q$, then (\ref{Poi}) (it is
called the {\it Poincar\'e inequality}\/ in this case) is valid for all domains such
that their boundaries are locally graphs of continuous functions in Cartesian
coordinates (see, for example, the classical book \cite{CH} by R.~Courant and D.~Hilbert
for the proof which can be easily extended from $p=2$ to any $p$).

Furthermore, if $p<n$ and $q=p^*$, then (\ref{Poi}) (it is called the {\it
Poincar\'e--Sobolev inequality}\/ in this case) holds for any bounded
$n$-dimensional Lipschitz domain. Moreover, the inequality is true provided $\Omega$ is a
{\it John domain}\footnote {This class of domains, more general than the Lipschitz ones, was
introduced by F.~John [22].} as was proved by B.~Bojarski \cite{Boj}.

Finally, we notice that if $q \neq p^*$, then (\ref{Poi}) holds if and only if
$L^{1,p} (\Omega)$ is continuously embedded into $L^{q}(\Omega)$. This was established by 
J.~Deny and J.-L.~Lions \cite{DL} for $p=q$; for the general case see \cite{NP1}.

Thus, the first point to be clarified about inequality (\ref{Poi}) concerns
smoothness of $\partial \Omega$. To a great extent, this is realised by V.\,G.~Maz'ya in
his comprehensive monograph {\it Sobolev Spaces}, where he presents his own results
and surveys those of other authors. (Recently, the 2nd revised and augmented edition
\cite{Maz} was published; its bibliography exceeds 800 entries. Moreover, several
sections deal with the question of exact constants in some inequalities.) Proofs of
basic facts can be found also in the recent textbook \cite{NP}; its English
translation is currently in preparation.\medskip

Almost everything known about sharp constants in various versions of inequality
(\ref{q-p}) comes under one of the following three conditions:
\begin{enumerate}
 \item $p=q=2$ (quadratic case);
 \item $\Omega=(0,l)$ (one-dimensional case);
 \item $p<n$, $q=p^*$ (critical case).
\end{enumerate}



It was mentioned above that the sharp constant in (\ref{q-p}) is $\lambda_1^{-1/2}$ in
the {\it quadratic case}. Here $\lambda_1 = \lambda^{\Omega}_1 \big( \lambda^{N}_1 \big)$
is the smallest positive eigenvalue of the Dirichlet (Neumann, respectively)
Laplacian for the Steklov (Poincar\'e, respectively) inequality. Explicit values of
these eigenvalues are found only for several particular domains. Among them, one
finds the following (see \cite{PS}).

\begin{center}
 \begin{tabular}{|l|c|c|}
 \hline \hfil Domain\vphantom{$\left(\frac {b^b}{b^b}\right)^2$}\hfil & $\lambda_1^D$ & $\lambda_1^N$ \\
 \hline Rectangle $a\times b$ \vphantom{$\left(\frac {b^b}{b^b}\right)^2$}
 & $\left( \frac{\pi}{a}\right)^2+\left( \frac{\pi}{b}\right)^2$ & $\big[ \frac{\pi}{\max\{a,b\}}\big]^2$\\
 \hline $45^\circ$ right triangle with the leg length $a$ \vphantom{$\left(\frac {b^b}{b^b}\right)^2$}
 & $5\left( \frac{\pi}{a}\right)^2$ & $\big( \frac{\pi}{a}\big)^2$\\
 \hline  $30^\circ$ right triangle with the hypotenuse length $a$ \vphantom{$\left(\frac {b^b}{b^b}\right)^2$}
 & $\frac{112}9\left( \frac{\pi}{a}\right)^2$ & $\frac{16}3\big( \frac{\pi}{a}\big)^2$\\
 \hline Equilateral triangle with the side length $a$ \vphantom{$\left(\frac {b^b}{b^b}\right)^2$}
 & $\frac{16}3\left( \frac{\pi}{a}\right)^2$ & $\frac{16}9\big( \frac{\pi}{a}\big)^2$\\
 \hline Disk of the radius $a$ \vphantom{$\left(\frac {b^b}{b^b}\right)^2$}
 & $\big( \frac{j_{0,1}}{a}\big)^2$ & $\big( \frac{j_{1,1}}{a}\big)^2$\\
 \hline
 \end{tabular}
\end{center}
%
%
%
%
%
Here $j_{0,1}$ $(j_{1,1})$ is the first positive zero of the Bessel function $J_0$
($J_1$, respectively). The Dirichlet and Neumann eigenvalues for sectors and annuli can also be
expressed in terms of Bessel functions.

Moreover, let $\Omega_1 \subset \RR^m$ and $\Omega_2 \subset \RR^n$, and let
$\lambda_1^{(j),D}$ and $\lambda_1^{(j),N}$ $(j=1,2)$ denote the fundamental
Dirichlet and Neumann, respectively, eigenvalues for the domains $\Omega_1$ and $\Omega_2$.
Then $\lambda_1^D = \lambda_1^{(1),D} + \lambda_1^{(2),D}$ and $\lambda_1^N = \min\{\lambda_1^{(1),N}, \lambda_1^{(2),N}\}$
are the corresponding eigenvalues for $\Omega = \Omega_1 \times \Omega_2$.\medskip

For estimating $\lambda_1^{D}$ one can use its monotonicity with respect to domain
variation and properties of the Steiner symmetrization (see \cite{PS}). In
particular, among all quadrilaterals of the same area the least value of
$\lambda_1^{D}$ is delivered by the square, whereas the equilateral triangle has the
least value of $\lambda_1^{D}$ among all triangles of the same area (see also
\cite{Fre}). Finally, a ball in $\RR^n$ has the least value of $\lambda_1^{D}$ among
all figures of the same area/volume\footnote{It must be emphasized that all estimates involving symmetrization
for their derivation are true for arbitrary $p$ and $q$. Thus, under the condition
that $u$ vanishes on $\partial \Omega$ the sharp constant in (\ref{q-p}) has the largest
value for a ball in $\RR^n$ (comparing other domains of the same area/volume).
Unfortunately, bounds for sharp constants are implicit unless $p=q=2$.}. In 1877, the two-dimensional version of the
last assertion was conjectured by Lord Rayleigh (see \cite{Ray}, pp.~339--340). It
was proved independently by G.~Faber \cite{Fa} and E.~Krahn \cite{Kr1}, \cite{Kr2}.

Less is known about estimates of the first positive Neumann eigenvalue. The
classical result of G.~Szeg\H o \cite{Sze} ($n=2$) and H.\,F.~Weinberger \cite{W}
(higher dimensions) says that a ball in $\RR^n$ has the largest value of
$\lambda_1^{N}$ among all domains of the same area/volume (see also \cite{AB}).
Analogous result for triangles was obtained recently in \cite{LS}. A global lower
bound for $\lambda_1^{N}$ was obtained for {\it convex} domains in \cite{PW};
namely, $\lambda_1^{N} > \big( \frac{\pi}{{\rm diam}\,\Omega} \big)^2$ unless $n=1$ when
$\Omega$ is an interval\footnote{A generalization of this result for eigenvalues of some nonlinear
Neumann problems was established recently in \cite{ENT}.}. There are also inequalities between the Dirichlet and
Neumann eigenvalues (see the recent paper \cite{Fi} for a brief historical survey).
Furthermore, it is shown in \cite{Sta} that if $p=q$ is arbitrary and (\ref{Poi})
holds for $\Omega_1 \subset \RR^m$ and $\Omega_2 \subset \RR^n$ with the sharp constants $C_1$
and $C_2$, respectively, then the sharp constant in the same inequality for $\Omega_1
\times \Omega_2$ is less than or equal to $\sqrt{2} \, (C_1+C_2)$.\medskip

We turn to the {\it one-dimensional} case and assume without loss of generality that
$\Omega=(0,1)$. If $u$ vanishes at the end-points, then the sharp constant in (\ref{q-p})
is as follows: 
\begin{equation}
C = C_1(p,q) = \frac{{\mathfrak F} \left( q^{-1} + p'^{-1} \right)}{2 \,{\mathfrak F}
\left( q^{-1} \right) {\mathfrak F}\left( p'^{-1} \right)} \, , \label{Schmidt}
\end{equation}
where ${\mathfrak F}(s) = \frac{\Gamma (s+1)}{s^s}$ and $p' = \frac{p}{p-1}$ is the
H\"older conjugate exponent to $p$. It was obtained by E.~Schmidt \cite{Sch} (the
case $p=q$ was considered earlier in \cite{Lev}; see also \cite{HLP}, sect.~7.6).
The extremal function, say, $U$ can be expressed in quadratures and is symmetric
with respect to $x - \frac{1}{2}$.\medskip

The one-dimensional Poincar\'e-type inequality has a more complicated story. It took
several years after the pioneering paper \cite{DGS}\footnote{In \cite{Maz}, 
sect.~1.1.19, the first result for $p=q$ is
attributed to A.~Stanoyevitch. However, the proof in his PhD thesis (1990) turned
out to be incorrect.} to establish the following
result (see \cite{BKN}, \cite{Naz} and also the recent paper \cite{GN} for a more
general problem and a historical survey):

\vspace{1.2mm}

\noindent {\it Let $n=1$ and $\Omega=(0, 1)$. If $q \leq 3 p$, then the sharp constant
in $(\ref{Poi})$ is equal to $C_1(p,q)$ defined by $(\ref{Schmidt}$ and the
corresponding extremal function $V$ is as follows: 
\begin{equation*}
 V(x)=U\big(x+\frac{1}{2}\big),\quad x\leq \frac{1}{2};\qquad  V(x)=-U\big(x-\frac{1}{2}\big),\quad x\geq \frac{1}{2},
\end{equation*}
where $U$ is the Schmidt function. In particular, $V$ is antisymmetric with respect
to $x - \frac{1}{2}$. Otherwise, the constant in $(\ref{Poi})$ is greater than
$C_1(p,q)$, and $V$ has no symmetry.}

\vspace{1.2mm}

Let us turn to the {\it critical}\/ case. In 1960, V.\,G.~Maz'ya \cite{Maz60} and 
H.~Federer and W.\,H.~Fleming \cite{FF} found the sharp constant in the Sobolev
inequality for $p=1$. It is equal to $\omega_{n-1}^{-1/n}\cdot n^{(1-n)/n}$, where
$\omega_{n-1} = 2 \, \pi^{n/2} / \Gamma (n/2)$ is the $(n-1)$-dimensional measure of
the unit sphere in $\RR^n$. Only fifteen years later, T.~Aubin \cite{Aub} and 
G.~Talenti \cite{Tal} obtained the exact constant for $p>1$ despite the fact that the
Bliss inequality \cite{Bl} and symmetrization --- the key ingredients of the
proof --- were known for a long time. It reads ($\cal B$ stands for Euler's
beta function): 
\begin{equation}
C = C_2 (n,p) = \omega_{n-1}^{-\frac{1}{n}} \, n^{-\frac{1}{p}} \Big(
\frac{p-1}{n-p} \Big)^{\frac{1}{p'}} \left[ {\cal B} \Big( \frac{n}{p} ,
\frac{n}{p'} + 1 \Big) \right]^{-\frac{1}{n}} , \label{Sob}
\end{equation}
and is {\it not attained} unless $\Omega = \RR^n$. In the paper \cite{CNV}, the constant
$C_2(n,p)$ was obtained by virtue of the mass transportation approach (the
generalized Monge--Kantorovich problem).\medskip

The situation is again more complicated for the Sobolev--Poincar\'e inequality. It
is known that for any John domain the sharp constant is greater than or equal to
$2^{1/n} C_2(n,p)$, where $C_2(n,p)$ is defined by (\ref{Sob}). Moreover, if $\Omega$ is
a ${\cal C}^2$-domain and $C$ in (\ref{Poi}) is strictly greater than $2^{1/n}
C_2(n,p)$, then the sharp constant {\it is attained} for this $\Omega$. In particular,
for any bounded ${\cal C}^2$-domain there exists $\beta>0$ such that the sharp
constant in the Sobolev--Poincar\'e inequality is attained when $1 < p <
\frac{n+1}{2} + \beta$ (see \cite{DN} for the proof). In the survey article
\cite{NazSurv}, the question when the sharp constant is attainable is discussed for
various critical inequalities.\medskip

In conclusion, we consider the following ``boundary analogue'' of inequality
(\ref{Poi}): 
\begin{equation}
\| u - \langle u \rangle_G \|_{L^q(G)} \leq C \|\nabla u \|_{L^p(\Omega)} , \qquad \langle u
\rangle_G = \frac{\int_G u \, \D S}{{\rm meas}_{n-1} \, G} . \label{PoinBd}
\end{equation}
Here $\Omega$ is a bounded Lipschitz domain in $\RR^n$, $n \geq 2$, $G$ is a part of
$\partial \Omega$ possibly coinciding with $\partial \Omega$.

Inequality (\ref{PoinBd}) holds for $u \in L^{1,p} (\Omega)$ provided 
\begin{equation*}
\aligned
&q\le p^{**}=\frac{(n-1)p}{n-p},&\quad\mbox{if}&\quad 1\le p<n;\\
&q<\infty,&\quad\mbox{if}&\quad p=n;\\
&q\le\infty,&\quad\mbox{if}&\quad p>n.
\endaligned
\end{equation*}

In the {\it quadratic} case (that is, $p=q=2$), the sharp constant in (\ref{PoinBd})
is again equal to $\lambda_1^{-1/2}$, but now $\lambda_1 = \lambda_1^S$ is the
smallest positive eigenvalue of the following mixed (unless $G = \partial \Omega$)
Steklov problem: 
\[ \Delta u = 0 \ \ \mbox{in} \ \Omega , \quad \frac{\partial u}{\partial {\bf n}} = \lambda u 
\ \ \mbox{on} \ G , \quad \frac{\partial u}{\partial {\bf n}} = 0 \ \ \mbox{on} \
\partial \Omega \setminus G .
\]
Here ${\bf n}$ is the exterior unit normal existing almost everywhere on $\partial
\Omega$. For $n = 2$ $(n = 3)$ and a particular choice of $G$, eigenvalues of the above
problem give {\it sloshing frequencies} of free oscillations of a liquid in a
channel (container, respectively); see, for example, \cite{Lamb}, ch. IX.

In \cite{NRep}, $\lambda_1^S$ is found for several simple domains with different
sets chosen as $G$. For example, let $\Omega$ be a $45^\circ$ right triangle with leg
equal to $a$, then:

\begin{center}
 \begin{tabular}{|l|c|}
 \hline \hfil $G$ \vphantom{$\left(\frac {b^b}{b^b}\right)^2$}\hfil  & $\lambda_1^S$ 
 \\
 \hline $G$ is the hypotenuse \vphantom{$\left(\frac {b^b}{b^b}\right)^2$} & $\frac{\sqrt{2}}{a}$ 
 \\
 \hline $G$ is a leg \vphantom{$\left(\frac {b^b}{b^b}\right)^2$} & $\frac{z_1^{(1)}\tanh(z_1^{(1)})}{a}\approx\frac{2.3236}{a}$ 
 \\
 \hline Two legs form $G$ \vphantom{$\left(\frac {b^b}{b^b}\right)^2$} & $\frac{2z_1^{(2)}\tanh(z_1^{(2)})}{a}\approx\frac{1.3765}{a}$ 
 \\
 \hline
 \end{tabular}
\end{center}
Here $z_1^{(1)}$ and $z_1^{(2)}$ are the smallest positive roots of the equations $\tan z + \tanh z = 0$ and $\tan z \, \tanh z = 1$, respectively.



In \cite{NRep} (see also \cite{Rep}), some applications of sharp constants
from (\ref{Poi}) and (\ref{PoinBd}) are considered. These applications concern
quantitative analysis of solutions and a posteriori error estimation for partial
differential equations.\medskip

In the {\it critical} case (that is, $p<n$, $q=p^{**}$), it should be emphasized
that the sharp constant in (\ref{PoinBd}) is related to that in the {\it trace
Sobolev inequality} for the half-space $\RR^n_+ = \{ x \in \RR^n\,:\ x_n > 0 \}$:
\begin{equation}
\| u (\cdot,0) \|_{L^{p^{**}} (\RR^{n-1})} \leq C_3(n,p) \| \nabla u
\|_{L^p(\RR^n_+)} , \label{trSob}
\end{equation}
valid for any $u \in L^{1,p} (\RR^n_+)$. In particular, it follows from
\cite{Maz88}, sect.~1.3, that $C_3(n,p) = 1$ for $p=1$.

J.\,F.~Escobar \cite{E} conjectured that if $p>1$ in (\ref{trSob}), then the
extremal function is equal to $|x-x^*|^{-(n-p) / (p-1)}$, where $x^* \notin \RR^n_+$
is arbitrary. This assertion is proved in \cite{E} only for $p=2$, but later the
general case was established in the remarkable paper \cite{Nt} based on the mass
transportation approach (see also \cite{Naz1}). This result implies that 
\[ C_3 (n,p) = \Big( \frac{p-1}{n-p} \Big)^{\frac{1}{p'}} \left[ 
\frac{\omega_{n-2}}{2} \, {\cal B} \Big( \frac{n-1}{2}, \frac{n-1}{2(p-1)} \Big)
\right]^{- \frac{1}{(n-1) p'}} .
\]

As in the case of the Sobolev--Poincar\'e inequality the following is true when the
case is critical. The sharp constant in (\ref{PoinBd}) is greater than or equal to
$C_3(n,p)$ for Lipschitz domains. Moreover, if $\Omega$ is a ${\cal C}^2$-domain and
$C>C_3(n,p)$, then the sharp constant {\it is attained} for this $\Omega$. In particular,
for any bounded ${\cal C}^2$-domain in $\RR^n$, $n \geq 3$, there exists $\delta > 0$
such that the sharp constant is attained for $1 < p < \frac{n+1}{2}+\delta$ (see
\cite{NRez} for the proof).

Since integral inequalities (as well as integration by parts) are at the heart of
theory of differential equations arising in mathematical physics, one might expect
that the interest to sharp constants in these inequalities will only intensify in
the future.

\end{document}